\documentclass[11pt]{article}

\usepackage{graphics}
\usepackage{amsfonts}
\usepackage{amssymb}
\usepackage{amsmath}
\usepackage{amsthm}
\usepackage{url}
\usepackage{color}
\usepackage{multirow}
\usepackage{enumerate}
\usepackage{hyperref}

\newtheorem{thm}{Theorem}

\newtheorem{cor}[thm]{Corollary}
\newtheorem{prop}[thm]{Proposition}

\theoremstyle{definition}

\theoremstyle{definition}

\theoremstyle{definition}

\newcommand{\R}{\mathbb{R}}

\DeclareMathOperator{\tr}{tr}
\DeclareMathOperator{\diag}{diag}
\DeclareMathOperator{\Diag}{Diag}

\title{Graph bisection revisited}
\author{R. Sotirov\thanks{Department of Econometrics and OR, Tilburg
University, The Netherlands. {\tt r.sotirov@uvt.nl} }}

\date{}

\begin{document}
\maketitle

\begin{abstract}
The graph bisection problem is the problem of partitioning the vertex set of a graph into two
sets of given sizes such that the sum of weights of edges joining  these two sets is optimized.
We present a   semidefinite programming  relaxation for the graph bisection problem with a matrix variable
of order $n$ - the number of vertices of the graph -
that is equivalent to the currently strongest semidefinite programming  relaxation obtained by using vector lifting.
The reduction in the size of the matrix variable enables us to impose additional
valid inequalities to the relaxation in order to further strengthen it.
The numerical results confirm that our simplified and  strengthened semidefinite relaxation
provides the currently strongest bound for the graph bisection problem in reasonable time.
\end{abstract}


\section{Introduction}\label{sec:intro}

The graph bisection  problem (GBP) is the problem of  dividing the vertices of a graph into two
sets of specified sizes such that the total weight of edges joining different sets is optimized.
The GBP is an NP-hard combinatorial optimization problem, see \cite{GarJoSt:76}.
It has many applications such as VLSI design \cite{Len:90}, parallel computing \cite{BisHenKa:00,HeKo:00,Simon91}, network partitioning
\cite{FidMatt:82,Sanch89}, and floor planing \cite{DaKuh87}. Graph partitioning also plays a role in machine learning (see e.g., \cite{LiSmola})
and data analysis (see e.g., \cite{{Pirim}}).

There are several SDP relaxations for the GBP  with  matrix variables of different orders.
In particular, there are relaxations whose matrices have orders $n$, $2n$, and $2n+1$, where $n$ is the order of the graph.
An SDP relaxation with a matrix variable of  order $n$ is introduced by Karisch, Rendl, and Clausen \cite{KaReCl:00}.
The same relaxation is  used by Feige and Langberg \cite{FeiLa:01}, and  Han, Ye, and Zhang \cite{HaYeZh:02} to derive approximation algorithms for the GBP.
Another SDP relaxation with a matrix variable of  order $n$ that is derived from an SDP relaxation for the more general graph partition problem
is introduced in \cite{Sot11}. In  \cite{Sot11} it is also proven that the above mentioned SDP relaxations of order $n$ are equivalent.

Wolkowicz and Zhao \cite{WolkZhao:99} derived an SDP relaxation with a matrix variable of  order $2n+1$.
This SDP relaxation with additional nonnegativity constraints dominates
the SDP relaxations with matrix variables of  order $n$, see \cite{KOP,Sot11}.

The GBP can be seen as a special case of the quadratic assignment problem (QAP).
De Klerk, Pasechnik, Sotirov, and Dobre  \cite{dKPaDoSo:10} exploited this to derive an SDP relaxation for the GBP from an SDP relaxation for the QAP,
which however reduces to a much smaller semidefinite program than the original QAP relaxation (see also \cite{KOP}).
In particular that relaxation contains matrix variables of orders $n$ and $2n$.
In \cite{Sot11}, it is proven  that the QAP-based SDP relaxation for the GBP is equivalent to the  strongest  SDP relaxation, that
is the SDP relaxation with nonnegativity constraints from \cite{WolkZhao:99}.

For specific families of (symmetric) graphs,
De Klerk et al.~\cite{dKPaDoSo:10} improved the QAP-based SDP relaxation for the GBP by adding a constraint that fixes one vertex of the graph.
Finally, in \cite{vDamSo:15} the SDP relaxation for the GBP  from \cite{dKPaDoSo:10}
was further strengthened by adding two constraints that correspond to assigning two vertices of the graph to different
parts of the partition. Both fixing-based  strengthening perform well on highly symmetric graphs.\\

In this paper, we present an SDP relaxation for the bisection problem whose matrix variable is of order $n$.
Our relaxation is equivalent to the strongest SDP relaxation for the GBP, that is the strongest Wolkowicz and Zhao relaxation  from \cite{WolkZhao:99}.
The new SDP relaxation exploits  the fact that the matrix variables corresponding to the two parts in the bisection are related.
Further,  we consider adding the facet defining inequalities of the boolean quadric polytope to our relaxation.
We also show that a large subset of the facet defining inequalities are redundant in the relaxation from  \cite{WolkZhao:99}.
The strengthened SDP bound outperforms all previously considered SDP bounds, including those tailored for highly symmetric graphs. \\

The paper is structured as follows. In Section \ref{sec:gpp} we provide an integer programming formulation of the problem, and
in Section \ref{ref:overview} an overview is given of the known SDP relaxations for the graph bisection problem.
In Section \ref{sect:revisites} we present our SDP relaxation and prove that it is equivalent to the strongest
SDP relaxation from \cite{WolkZhao:99}.
We further suggest how to improve our relaxation. Finally, in Section \ref{sect:numerres} we present  numerical results.

\section{The graph bisection problem}\label{sec:gpp}

In this section we formulate the minimum   graph bisection problem as an integer optimization problem.
Let $G=(V,E)$ be an undirected graph  with vertex set $V$, where $|V|=n$ and edge set $E$.
The goal is to find a partition  of the vertex set into two disjoint subsets $S_1$ and $S_2$ of
specified sizes $m_1\geq  m_2$, $m_1+m_2=n$ such that the sum of weights of edges joining
$S_1$ and $S_2$ is minimized. If  $m_1=m_2$  then one refers to the associated problem  as the  graph equipartition problem.
We consider here only the case that $m_1 > m_2$. For detailed analysis of the SDP relaxations for the   graph equipartition problem,  see \cite{Sot10}.

Let us denote by $A$ the adjacency matrix of $G$. For a given partition of the graph $G$ into two subsets,   let   $Z=(z_{ij})$ be
the $n\times 2$ matrix defined by
\[
z_{ij} := \left \{
\begin{array}{ll}
1 &  \mbox{if $i\in S_j$ } \\
0 & \mbox{otherwise}
\end{array} \right . \quad i=1,\ldots, n,~j=1,2.
\]
The $j$th column of $Z$ is the characteristic vector of $S_j$.
The  cut of the partition, which is the sum of weights of edges joining different sets,
 is equal to:
\[
\frac{1}{2} \tr A(J -ZZ^{\mathrm{T}}) = \frac{1}{2} \tr (L ZZ^{\mathrm{T}}),
\]
where  $L=\Diag(Ae) -A$ is the Laplacian matrix of the graph, and $J$ (resp.~$e$) the all-ones matrix (resp.~vector).
Therefore, the minimum  GBP problem can be formulated as follows
\begin{equation} \label{GPP}
\min  \left \{  \frac{1}{2} \tr (LZZ^{\mathrm{T}}):  Ze =e, ~Z^{\mathrm{T}} e =m, ~ z_{ij}\in \{0,1\},~~\forall i,j \right \},
\end{equation}
where $m=(m_1, m_2)^\mathrm{T}$.

\section{Overview of SDP relaxations} \label{ref:overview}

In this section we provide an overview of existing SDP relaxations for the GBP.
The following SDP relaxation is  derived in \cite{Sot11}
\begin{equation} \label{eq:RS}
\begin{array}{rl}
\min & \frac{1}{2} \tr (LX) \\[1ex]
{\rm s.t.} & \diag(X)= e, ~ \tr (JX) = m_1^2 +m_2^2  \\[1ex]
& 2X - J\succeq 0, ~X\in {\mathcal S}_n,
\end{array}
\end{equation}
where the `diag' operator maps an $n\times n$ matrix to the $n$-vector given by its diagonal,
and ${\mathcal S}_n$ denotes the space of $n\times n$ symmetric matrices.
Nonnegativity constraints on the matrix variable  in \eqref{eq:RS} are redundant.
This follows from $\diag(X)= e$ and $2X - J\succeq 0$, see \cite{vDamSo:15} for details.
The SDP relaxation \eqref{eq:RS} is equivalent to the SDP relaxation with a matrix variable of order $n$ from \cite{KaReCl:00}.

The following SDP relaxation for the GBP is derived in \cite{WolkZhao:99}:
\begin{equation}\label{eq:ZW}
\begin{array}{rl}
\min &  \frac{1}{2} \tr L(Y_{11}+Y_{22}) \\[1ex]
{\rm s.t.} &  \tr (Y_{ii}) = m_{i}, ~\tr (J Y_{ii}) = m_{i}^{2},  ~~i=1,2 \\[1ex]
&\diag(Y_{12})=0, ~\tr J(Y_{12}+Y_{12}^{\mathrm T})= 2m_{1}m_{2} \\[1ex]
&  Y = \left(
\begin{array}{cc}
Y_{11} & Y_{12}  \\
Y_{12}^{\mathrm T} & Y_{22} \\
\end{array} \right)
~~y=\diag(Y),
~~Y-yy^{\mathrm T} \succeq 0,~~ Y\geq 0,
\end{array}
\end{equation}
where $Y\in {\mathcal S}_{2n}$. From now on, we assume that matrices of order $2n$ have the block structure as given above.

Although the nonnegativity constraints were not included in the relaxation from \cite{WolkZhao:99},
the authors  mentioned that it would be worth adding them.
The SDP relaxation \eqref{eq:ZW} does not have strictly feasible solutions.
From a computational point of view, this is an indication that the model may be difficult to solve directly as it is.
Therefore,  Wolkowicz and Zhao  \cite{WolkZhao:99} derive the Slater feasible version of the relaxation whose matrix variable is of order $n$.
However,  that model includes multiplications with projection matrices of size $(2n+1) \times n$.
The above relaxation can be further strengthened by adding the following inequalities
\begin{eqnarray}
0 \leq y_{i,j} \leq y_{i,i} \label{BQP1} \\[1ex]
y_{i,i} + y_{j,j} \leq 1 + y_{i,j} \label{BQP2} \\[1ex]
y_{i,k} + y_{j,k} \leq y_{k,k} + y_{i,j} \label{BQP3} \\[1ex]
y_{i,i} + y_{j,j} + y_{k,k} \leq y_{i,j}+y_{i,k} + y_{j,k} +1, \label{BQP4} 
\end{eqnarray}
where $Y=(y_{ij})$ and  $ 1\leq  i,j,k \leq 2n$, $i\neq j$, $i\neq k$, $j\neq k$.
The  inequalities  \eqref{BQP1}--\eqref{BQP4}  are facet defining
inequalities of the boolean quadric polytope (BQP), see \cite{Pad:89}.

Wolkowicz and Zhao \cite{WolkZhao:99}  prove that
for a matrix $Y$  that is feasible for  the SDP relaxation \eqref{eq:ZW} the following is satisfied:
\begin{equation} \label{property}
Y_{11}+Y_{12}=y_1 e^{\mathrm T}, \quad Y_{12}^{\mathrm T} + Y_{22} = y_2 e^{\mathrm T}, \quad y_1+y_2= e, \quad Y_{ii}e=m_i y_i ~(i=1,2),
\end{equation}
where $y_i=\diag(Y_{ii})$ ($i=1,2$).
From here it follows that for given $Y_{11}$ and $y_1$ the above equations uniquely determine $Y_{12}, Y_{22}$ and $y_2$.
We will exploit  this  to derive the simplified SDP relaxation in the following section.

Extensive numerical results in \cite{Sot11} show that \eqref{eq:ZW} provides the strongest  SDP relaxation for the GBP.
To the best of our knowledge, we are not aware of numerical test that involve the SDP relaxation \eqref{eq:ZW} and the inequalities  \eqref{BQP1}--\eqref{BQP4}.

Finally, we prove that the optimal value of  the SDP relaxation  \eqref{eq:ZW} is at least that of the relaxation \eqref{eq:RS}.
\begin{prop}
Let $m_1>m_2$ and $m_1+m_2=n$. Then the SDP relaxation  \eqref{eq:ZW} dominates the SDP relaxation  \eqref{eq:RS}.
\end{prop}

\noindent
{\em Proof.} Let $Y_{ij}$ and $y_i$ ($i,j=1,2$) be feasible for  \eqref{eq:ZW}, and set $X =Y_{11} + Y_{22}$.
Now, $\tr (JX) = m_1^2 +m_2^2$ and $\diag(X)= e$ follow from feasibility of $Y_{ii}$ ($i=1,2$) and  \eqref{property}.
The SDP constraint follows from $\begin{pmatrix} Y_{ii} & y_i \\ y_i^{\mathrm T} & 1 \end{pmatrix}  \succeq 0$, $i=1,2$.
\qed

The similar result is proven in \cite{KOP}.
In particular,   it was proven that the QAP-based SDP relaxation for the GBP dominates the SDP relaxation from \cite{KaReCl:00}.
However, the QAP-based SDP relaxation for the GBP  is equivalent to  \eqref{eq:ZW},
and  the relaxation from \cite{KaReCl:00}  to \eqref{eq:RS}, see \cite{Sot11}.

\section{A simplified SDP relaxation } \label{sect:revisites}

In this section we derive an SDP relaxation for the GBP with a matrix variable of order $n$,
and prove that it is equivalent to the best known SDP relaxation for general graphs that is derived in \cite{WolkZhao:99}.
To derive the relaxation we exploit the fact that  the variables associated to the two sets in the bisection are related.
Namely, variables coming from the assignment to the second set are redundant in the assignment constraints.
It is surprising that this observation was not earlier exploited in the context of the GBP.
However, similar idea was used  in \cite{ReSo:16} to derive an SDP relaxation  for the vertex separator problem.

Our observation lead us to the following SDP relaxation:
\begin{equation} \label{newN}
\begin{array}{rl}
\min & \tr L( 2X +J - x e^{\mathrm T} - e x^{\mathrm T}) \\[1.5ex]
{\rm s.t.} & x^{\mathrm T} e = m_1, ~~ \tr(J X)= m_1^2, ~~X e = m_1 x ~~\\[1.5ex]
& X \geq 0, ~~  x e^{\mathrm T} - X \geq 0, ~~ J + X - x e^{\mathrm T} - e x^{\mathrm T} \geq 0 \\[1.5ex]
& X  \succeq 0, ~~\diag(X)=x, ~~X\in {\mathcal S}_n.
\end{array}
\end{equation}
All equality constraints in \eqref{newN} are related to the variables associated to the set $S_1$.
The constraints $X \geq 0$ ensure that the matrix variable corresponding to  $S_1$ is nonnegative, while
constraints $x e^{\mathrm T} - X \geq 0$, $J + X - x e^{\mathrm T} - e x^{\mathrm T} \geq 0$ do the same for the slack matrix variables.

One may wish to replace the  semidefinite constraint  $X  \succeq 0$ from \eqref{newN} with
the  in general stronger  constraint $ X -xx^{\mathrm T} \succeq 0$.
However, from the following result it follows that in our case those two semidefinite constraints are  equivalent.
\begin{prop} {\rm (\cite{Gij05}, Proposition 7)} \label{Gijs}
Let $X$ be a symmetric matrix of order $n$ such that $c\diag(X)=Xe$ for some $c\in \R$, and
\[
\bar{X}=\left (
\begin{array}{cc} 1 &~~ \diag(X)^{\mathrm{T}}\\[1ex]
\diag(X) & X  \end{array} \right ).
\]
Then the following are equivalent: \\[1ex]
\begin{tabular}{ll}
(i)& $\bar{X}$ is positive semidefinite, \\
(ii)& $X$ is positive semidefinite and $\tr (JX)\geq (\tr X)^2$.
\end{tabular}
\end{prop}
The equivalence of the two SDP constraints follows from the fact that for a feasible $X$ for \eqref{newN} one has $\tr(J X)=(\tr X)^2=m_1^2$ and $Xe=m_1 \diag(X)$.
We prove now our main result.
\begin{thm} \label{Thm:main}
Let $m_1+m_2=n$, $m_1>m_2$.
The SDP relaxations \eqref{eq:ZW}  and  \eqref{newN} are equivalent.
\end{thm}

\noindent
{\em Proof.} Let $X$ be feasible for  \eqref{newN} and $x=\diag(X)$.
We construct a feasible  $Y$, $y=\diag(Y)$ for \eqref{eq:ZW}  in the following way.
Define  $y_1:=x$, $y_2:= e-x$, $y:=(y_1,y_2)^{\mathrm T}$,  matrices
\[
Y_{11} := X, ~~~Y_{22} := J+X - x e^{\mathrm T} - e x^{\mathrm T}, ~~    Y_{12} := x e^{\mathrm T} - X,
\]
and collect all blocks into the matrix
\[
Y =
\left(
\begin{array}{cc}
Y_{11} & Y_{12}  \\
Y_{12}^{\mathrm T} & Y_{22} \\
\end{array} \right)
=\left(
\begin{array}{cc}
X &  x e^{\mathrm T} - X   \\
e x^{\mathrm T} - X & J+X - x e^{\mathrm T} - e x^{\mathrm T} \\
\end{array} \right ).
\]
Now, we first prove that
\[
\left(
\begin{array}{cc}
X &  x e^{\mathrm T} - X   \\
e x^{\mathrm T} - X & J+X - x e^{\mathrm T} - e x^{\mathrm T} \\
\end{array} \right )
-
\left(
\begin{array}{cc}
 x x^{\mathrm T} & x (e -x)^{\mathrm T} \\
 (e -x)x^{\mathrm T}  &  (e -x)(e -x)^{\mathrm T}
\end{array} \right )
\succeq 0.
\]
To show this, we rewrite  the left hand side of the matrix inequality above  as it follows
\[
\left(
\begin{array}{cc}
 X - x x^{\mathrm T}  &  x x^{\mathrm T} - X   \\
x x^{\mathrm T} - X  &  X - x x^{\mathrm T} \\
\end{array} \right ).
\]
Now, for arbitrary vectors $z_1, z_2 \in \R^n$ we have
\[
(z_1^{\mathrm T},z_2^{\mathrm T} )\left(
\begin{array}{cc}
 X - x x^{\mathrm T}  &  x x^{\mathrm T} - X   \\
x x^{\mathrm T} - X  &  X - x x^{\mathrm T} \\
\end{array} \right )
\left (
\begin{array}{c}
z_1 \\ z_2
\end{array}
\right ) = (z_1-z_2)^{\mathrm T} (X-xx^{\mathrm T} )(z_1-z_2) \geq 0,
\]
from where it follows the claim.
Let us now verify $\tr(JY_{22})=m_2^2$. Namely,
\[
\tr(JY_{22}) = \tr(J(J+X - x e^{\mathrm T} - e x^{\mathrm T}))=n^2+m_1^2-2nm_1=m_2^2.
\]
Similarly, the remaining constraints from \eqref{eq:ZW}  can be  verified.

Conversely, let $Y$ be feasible for  \eqref{eq:ZW}. We set $X = Y_{11}$ and $x =\diag(Y_{11})$.
Since every feasible matrix  $Y\in {\mathcal S}_{2n}$ for  \eqref{eq:ZW} satisfies also \eqref{property}, feasibility of $X$ follows by direct verification.
Finally, it is not difficult to check that the objectives coincide for any pair of feasible solutions $(Y,X)$. \qed\\

Note that the result from the previous theorem is also valid when $m_1=m_2$.
However, it was proven in \cite{Sot10} that all known vector and matrix lifting based SDP relaxations
for the  $k$-equipartition problem  ($k\geq 2$)  are equivalent.

It is not difficult to verify that the SDP relaxation  \eqref{newN} has a strictly feasible point. In deed, the following matrix is in the interior of the
feasible set of \eqref{newN}:
\[
\hat{X} = \frac{m_1}{n}I + \frac{m_1(m_1-1)}{n(n-1)}(J-I),
\]
where $I$ is the identity matrix.

The following result is a direct consequence of Theorem \ref{Thm:main}.
\begin{cor}\label{cor:Old}
The SDP relaxation \eqref{eq:ZW} without nonnegativity constraints
 is equivalent to the SDP relaxation  \eqref{newN} without nonnegativity constraints, i.e.,
\[
\begin{array}{rl}
\min & \tr L( 2X +J - x e^{\mathrm T} - e x^{\mathrm T}) \\[1.5ex]
{\rm s.t.} & x^{\mathrm T} e = m_1, ~~ \tr(J X)= m_1^2, ~~X e = m_1 x ~~\\[1.5ex]
& X \succeq 0, ~~\diag(X)=x, ~~X\in {\mathcal S}_n.
\end{array}
\]
\end{cor}
In order to improve the SDP relaxation  \eqref{newN} we can add  the facet defining
inequalities of the boolean quadric polytope,  see \cite{Pad:89}.
We first note that the inequality constraints
$X \geq 0,$ $x e^{\mathrm T} - X \geq 0,$ and $J + X - x e^{\mathrm T} - e x^{\mathrm T} \geq 0$
from  the SDP relaxation \eqref{newN} are exactly  the following BQP constraints
\[
0 \leq x_{i,j} \leq x_{i,i}, ~~x_{i,i} + x_{j,j} \leq 1 + x_{i,j}, ~~ i,j=1,\ldots, n, ~i\neq j.
\]
Note also that the SDP relaxation from Corollary \ref{cor:Old} differs from  the SDP relaxation \eqref{newN} exactly for those constraints.
Thus, in order to strengthen the SDP relaxation \eqref{newN} one can add the following BQP constraints:
\begin{equation}\label{cutN}
x_{i,k} + x_{j,k} \leq x_{k,k} + x_{i,j}, \quad
x_{i,i} + x_{j,j} + x_{k,k} \leq x_{i,j}+ x_{i,k} + x_{j,k} +1,
\end{equation}
for $ 1\leq  i,j,k \leq n$, $i\neq j$, $i\neq k$, $j\neq k$.\\

Let us now show that the bound obtained by solving the SDP relaxation \eqref{eq:ZW} with additional BQP constraints \eqref{BQP1}--\eqref{BQP4}
is equal to the bound obtained by solving \eqref{newN}  with additional constraints \eqref{cutN}.
Let $Y=(y_{ij})$ ($i,j=1,\ldots, 2n$) be feasible for  \eqref{eq:ZW}.
To show that the following inequalities
\[
y_{i,i} + y_{j,j} \leq 1 + y_{i,j}, ~~n+1\leq  i,j \leq 2n, ~i\neq j,
\]
are redundant,  it is instructive to look at \eqref{property}.
From \eqref{property} we have that $Y_{22}=J+Y_{11} - y_1e^{\mathrm T} - y_1 e^{\mathrm T}$,
and therefore  the above inequalities reduce to the redundant constraints $y_{ij}\geq 0$ ($1\leq i,j \leq n$).

To show that
\[
 y_{i,k} + y_{j,k} \leq y_{k,k} + y_{i,j}, ~~ 1\leq i \leq n, ~n+1\leq  j,k \leq 2n,~ j\neq k,
\]
we again consider \eqref{property} and obtain
\[
y_{i,j}+y_{k,j}\leq y_{j,j}+y_{i,k},~~ 1\leq i,j,k \leq n.
\]
In a similar way we get that the only non-redundant constraints among \eqref{BQP1}--\eqref{BQP4} are exactly of the form \eqref{cutN}.
This is summarized as follows.
\begin{thm}
The SDP relaxation  \eqref{newN}  with additional constraints \eqref{cutN} is equivalent to the SDP relaxation
\eqref{eq:ZW} with additional BQP constraints \eqref{BQP1}--\eqref{BQP4}.
\end{thm}

This paper does not only present reformulated and simplified the currently strongest SDP relaxation for the the bisection problem,
but also suggests its further strengthening.
In the following section we test our simplified and strengthened SDP relaxation on several graphs from the literature.

\section{Numerical results} \label{sect:numerres}

In this section we present numerical results that verify the  quality of the SDP relaxation \eqref{newN},
as well as the relaxation obtained after adding the BQP constraints \eqref{cutN} to \eqref{newN}.
All relaxations were solved with Mosek \cite{mosek} using the Yalmip interface \cite{YALMIP}
on an  Intel Xeon, E5-1620, 3.70 GHz  with 32 GB memory.

The instances we use belong to the various classes of graphs from the literature.
In particular, in Table \ref{tab1} we consider the following graphs.
\begin{itemize}
\item {\tt compiler design instances} were introduced in \cite{JoMeNe93}.
We denote them by  {\tt cd.xx.yy}.

\item {\tt kkt instances}  originate from nested dissection  approaches for solving sparse symmetric linear systems,  see \cite{Helm04}.
We denote them by {\tt kkt}\_{\tt name}.

\item {\tt mesh instances} come from an application of the finite element methods, see \cite{SoKeWoZo}.
We denote them with the initials {\tt mesh.xx.yy}.

\item {\tt VLSI design instances} are derived  from data in the layout of electronic circuits. For details see
\cite{FeMaSoWeWo:98}. We denote them with the initials {\tt vlsi.xx.yy}.
\end{itemize}
In the above instances {\tt xx} denotes the number of vertices, and {\tt yy} the number of edges in the graph.
Table \ref{tab1} reads as follows.
In the first  three columns, we list the graphs,  number of vertices in the graph, and corresponding $m$, respectively.
In the fourth to six column we present the SDP bounds \eqref{eq:RS}, \eqref{newN}, and the SDP bound \eqref{newN}
with additional BQP constraints \eqref{cutN}, respectively.
Bounds in the column six are  obtained by adding the most violated inequalities of type \eqref{cutN} to the SDP relaxation \eqref{newN}.
The cutting plane scheme adds at most $2n$ violated valid constraints in each iteration and performs at most 20 iterations.
In the last column of Table \ref{tab1}  we list upper bounds obtained by a tabu search heuristics.

All lower bounds  in Table \ref{tab1}  are rounded up to the closest  integer.
Note that for only three out of twenty-one instances  we can not prove optimality.
We compute the bound \eqref{newN} for ${\tt kkt}\_{\tt putt}01 $ ($n=115$)   in 106 seconds, and  \eqref{newN}+\eqref{cutN} in 21 minutes.
To prove optimality for e.g., ${\tt mesh}.70.120 $ we need less than one minute.
\\
\begin{table}[h!]
\begin{center}
{\begin{tabular}{|l|c|c|ccc|c|}
\hline
instance & $|V|$ &  $m^\mathrm{T}$  & \eqref{eq:RS}  & \eqref{newN} & \eqref{newN}+\eqref{cutN} & u.b. \\[1ex]
\hline
{\tt cd}.30.47 & 30 &  $(20,10)$  &  110 & 114 &  114 & 114 \\
{\tt cd}.30.56 & 30  &  $(20,10)$  & 156 & 169 & 169 & 169     \\
{\tt cd}.45.98 & 45 &  $(25,20)$  & 576 & 631 & 631 & 631 \\
{\tt cd}.47.99 & 47 &  $(25,22)$ &  471 & 514 & 537 & 537  \\
{\tt cd}.47.101 & 47 & $(25,22)$ & 326 & 361 & 382 & 382    \\
{\tt cd}.61.187 & 61 & $(40,21)$ &   774 & 798 &  798 &  798  \\ \hline

{\tt kkt}\_{\tt lowt}01 & 82 & $(42,40)$ & 5 &  5 &  13  & 13 \\
{\tt kkt}\_{\tt putt}01 & 115 & $(59,56)$ & 20 & 22 &  28 &  29 \\ \hline

{\tt mesh}.35.54 & 35 & $(22,13)$ & 2&  4 & 4   & 4 \\
{\tt mesh}.69.212 & 69 & $(40,29)$ & 2   & 2  & 4 & 4 \\
{\tt mesh}.70.120 & 70 & $(50,20)$ & 2  & 4 & 6 & 6  \\
{\tt mesh}.74.129 & 74 & $(70,4)$ & 1 & 4   & 4 &4  \\
{\tt mesh}.137.231 & 137 & $(100,37)$ & 1& 3& 6 & 6   \\
{\tt mesh}.148.265 & 148 & $(120,28)$ & 1 & 5 &  6 & 6  \\  \hline
%
{\tt vlsi}.34.71 & 34 & $(22,12)$ & 4 & 6   &6  & 6  \\
{\tt vlsi}.37.92 & 37 & $(30,7)$ & 3 & 6   & 6 &  6 \\
{\tt vlsi}.38.105 & 38 & $(20,18)$& 84  & 86  &110  &110 \\
{\tt vlsi}.42.132 & 42 & $(20,22)$& 97 & 99  & 120   & 120  \\
{\tt vlsi}.48.81 & 48 & $(40,8)$ & 4 & 12  & 12 & 18  \\
{\tt vlsi}.166.504 & 166 & $(100,66)$ & 12 & 23 & 24 & 24 \\
{\tt vlsi}.170.424 & 170 & $(100,70)$& 35 & 37 &  37& 48  \\\hline
\end{tabular}}
\caption{Computational results for the bisection problem. }\label{tab1}
\end{center}
\end{table}

In \cite{vDamSo:15}, the authors  strengthened the SDP relaxations  \eqref{eq:RS} and \eqref{eq:ZW} by adding
two constraints that correspond to assigning two vertices of the graph to different parts of the partition.
In particular, they  show that such strengthening performs well on highly symmetric graphs
when other relaxations provide weak or trivial bounds.
In \cite{vDamSo:15}, it  was also shown how to aggregate
the triangle and independent set constraints for highly symmetric graphs in order to add them to the SDP relaxation \eqref{eq:RS}.
Our numerical results show that the SDP relaxation \eqref{newN} with additional inequalities \eqref{cutN}
provides bounds that are competitive to those from \cite{vDamSo:15}.

In particular, in Table \ref{tab2} we list bounds for highly symmetric graphs considered in \cite{vDamSo:15}.
{\tt Pappus}, {\tt Desargues},  and {\tt Biggs-Smith} graphs  are distance-regular graphs,  $J(7,2)$ is the Johnson graph.
The first three columns in Table \ref{tab2} read similar as the first three columns in Table \ref{tab1}.
In the fourth (resp.~sixth) column we list  values of the SDP bound \eqref{newN} (resp.~bound \eqref{newN} with additional
  inequalities \eqref{cutN}) for different graphs.
The fifth column of Table \ref{tab2} lists the best obtained bounds from \cite{vDamSo:15};
that is for the  {\tt Pappus} graph the relaxation \eqref{eq:RS} with all triangle inequalities,
for {\tt Desargues}    the relaxation  \eqref{eq:ZW} with constraints that fix two vertices of the graph,
for $J(7,2)$  the relaxation  \eqref{eq:RS} with independent set inequalities,
and for {\tt Biggs-Smith} the relaxation  \eqref{eq:RS} with all triangle inequalities.
To compute the SDP bound \eqref{newN}  for   highly symmetric graphs
we didn't exploit symmetry reduction as described in \cite{vDamSo:15} although this can be done in a similar way.
By doing as described in \cite{vDamSo:15}, one can compute bounds from Table \ref{tab2} very fast. The interested reader is invited to verify this.
\begin{table}[h!]
\begin{center}
\begin{tabular}{|c|c|c|ccc|c|}
\hline
$G$ & $|V|$ & $m^\mathrm{T}$ & \eqref{newN} & b.b.~\cite{vDamSo:15}  & \eqref{newN}+\eqref{cutN}  & u.b. \\[1ex]
\hline
{\tt Pappus} & 18   & $(10,8)$  & 6  & 7 & 7  & 8   \\

{\tt Desargues}& 20 & $(15,5)$  & 5  & 6 & 6  & 7  \\

$J(7,2)$      &21  & $(11,10)$ & 37 & 40 & 40 & 40  \\

{\tt Biggs-Smith}& 102& $(70,32)$ & 10 & 15 & 15 & 18 \\\hline
\end{tabular}
\caption{Bounds for the bisection on highly symmetric graphs. }  \label{tab2}
\end{center}
\end{table}

\section{Conclusion}
In this paper we present an SDP relaxation for the graph bisection problem with a matrix variable of order $n$, where $n$ is the order of the graph.
To derive our relaxation we exploit the fact that variables corresponding to one set in the bisection uniquely determine variables
of the other set. We prove that our relaxation is equivalent to the strongest known SDP relaxation for general graphs
that is obtained by using vector lifting.
This result is in the line of the similar results for some other optimization problems.
Namely,   for the graph equipartition problem there exists an SDP relaxation with a matrix variable of order equal to
the order of the graph, which is equivalent to the strongest vector lifting-based SDP relaxation, see \cite{Sot10}.

To strengthen our SDP relaxation we add facet defining inequalities of the boolean quadric polytope, which
enables us to compute strongest SDP bounds for the GBP and for graphs with $n\leq 200$ vertices in reasonable time.

Since our relaxation has strictly feasible solutions it can be directly solved as it is,
which makes it attractive for a branch and bound framework.
However, this will be part of our future research.


\begin{thebibliography}{10}

\bibitem{BisHenKa:00}
R. Biswas, B. Hendrickson, G. Karypis.
Graph partitioning and parallel computing. {\em Parallel Comput.},  26(12):1515--1517, 2000.

\bibitem{DaKuh87}
W. Dai, E. Kuh. Simultaneous floor planning and global routing for hierarchical building-block layout.
 IEEE Trans. Comput.-Aided Des. {\em Integrated Circuits $\&$ Syst.}, CAD-6, 5:828--837, 1987.

\bibitem{vDamSo:15}
E.R. van Dam, R. Sotirov. Semidefinite programming and eigenvalue bounds for the graph partition problem,
{\em Math. Program., Ser B}, 151(2):379--404, 2015.

\bibitem{dKPaDoSo:10}
E.\ De Klerk, D.V.\ Pasechnik, R.\ Sotirov,  C.\ Dobre.
On semidefinite programming relaxations of maximum $k$-section. {\em Math.  Program. Ser. B}, 136(2):253--278, 2012.

\bibitem{KOP}
E.~de Klerk, F.M.~de Oliveira Filho,  D.V.~Pasechnik.
 Relaxations of combinatorial problems via association   schemes, in
 {\em Handbook of Semidefinite, Cone and Polynomial Optimization: Theory, Algorithms, Software and Applications},
  Miguel Anjos and Jean Lasserre (eds.), pp. 171--200, Springer, 2012.

\bibitem{FeiLa:01}
U. Feige,  M. Langberg.
Approximation algorithms for maximization problems arising in graph partitioning. {\em J. Algorithm}, 41:174--211, 2001.

\bibitem{FeMaSoWeWo:98}
C.E. Ferreira, A. Martin, C.C. de Souza, R. Weismantel, L.A. Wolsey.
The node capacitated graph partitioning problem: a computational study. {\em Math. Program.}, 81:229--256, 1998.

\bibitem{FidMatt:82}
C.M. Fiduccia,  R.M. Mattheyses. A linear-time heuristic for improving network partitions.
Proceedings of the 19th Design Automation Conference, 175--181, 1982.

\bibitem{GarJoSt:76}
M.R. Garey,  D.S. Johnson, L. Stockmeyer. Some simplified NP-complete graph problems. {\em Theoret. Comput. Sci.}  1(3):237–-267, 1976.

\bibitem{Gij05}
D. Gijswijt. Matrix algebras and semidefinite programming techniques for codes.
PhD thesis, University of Amsterdam, The Netherlands, 2005.

\bibitem{HaYeZh:02}
Q. Han, Y. Ye, J. Zhang. An improved rounding method and semidefinite relaxation for graph partitioning. {\em Math. Program.}, 92:509--535, 2002.

\bibitem{Helm04}
C. Helmberg. A cutting plane algorithm for large scale semidefinite relaxations.
In: {\em Padberg Festschrift "The Sharpest Cut"}, M. Gr\"otschel (Ed.).  MPS-SIAM 233--256, 2004.

\bibitem{HeKo:00}
B. Hendrickson,  T.G. Kolda. Partitioning rectangular and structurally nonsymmetric sparse matrices for parallel processing.
 {\em SIAM J. Sci. Comput.}, 21(6):2048--2072, 2000.

\bibitem{JoMeNe93}
E. Johnson, A. Mehrotra, G. Nemhauser. Min-cut clustering. {\em Math. Program.}, 62:133--152,1993.

\bibitem{KaReCl:00}
S.E. Karisch, F. Rendl,  J. Clausen. Solving graph bisection problems with semidefinite programming.
{\em INFORMS J. Comput.},  12:177--191, 2000.

\bibitem{Len:90}
T. Lengauer. {\em Combinatorial algorithms for integrated circuit layout}. Wiley, Chicester, 1990.

\bibitem{LiSmola}
M. Li, D.G. Andersen, A.J. Smola.
Graph partitioning via parallel submodular approximation to accelerate distributed machine learning, {\tt  arXiv:1505.04636v1}, 2015.

\bibitem{YALMIP}
J. L\"ofberg.
YALMIP: A toolbox for modeling and optimization in MATLAB.   {\it Proc. of the CACSD Conference, Taipei, Taiwan}, 284--289, 2004.
\url{http://users.isy.liu.se/johanl/yalmip/}

\bibitem{mosek}
 MOSEK ApS. The MOSEK optimization toolbox for MATLAB manual. Version 7.1 (Revision 28), 2015.
 \url{http://docs.mosek.com/7.1/toolbox/index.html}

\bibitem{Pad:89}
 M.W. Padberg.  The boolean quadric polytope: Some characteristics, facets and relatives. {\em Math. Program.},  45:139--172, 1989.

\bibitem{Pirim}
H. Pirim, B.  Ek\c{c}io\u{g}lu,  A. Perkins,  \c{C}etin Y\"uceer.
Clustering of high throughput gene expression data. {\em Comput Oper Res.}, 39(12):3046-–3061, 2012.

\bibitem{ReSo:16}
F. Rendl, R. Sotirov. The min-cut and vertex separator problem, Preprint 2016.

\bibitem{Sanch89}
L. Sanchis. Multiple-way network partitioning. {\em IEEE Trans. Comput.},  38:62--81, 1989.

\bibitem{Simon91}
H.D. Simon. Partitioning of unstructured problems for parallel processing.  {\em Comput. Syst. Eng.},  2:35--148, 1991.

\bibitem{Sot10}
R. Sotirov. SDP relaxations for some combinatorial optimization problems.
In M.F. Anjos and J.B. Lasserre (Eds.),
{\em Handbook of Semidefinite, Conic and Polynomial Optimization: Theory, Algorithms, Software and Applications}, pp.~795-820, 2012.

\bibitem{Sot11}
R. Sotirov. An efficient semidefinite programming relaxation for the graph partition problem, {\em INFORMS J. Comput.}, 26(1):16--30, 2014.

\bibitem{SoKeWoZo}
C.C. de Souza, R. Keunings, L.A. Wolsey, O. Zone. A new approach to minimizing the frontwidth in finite element calculations.
{\em  Computer Methods in Applied Mechanics and Engineering}, 111:323–334, 1994.

\bibitem{WolkZhao:99}
H. Wolkowicz, Q. Zhao. Semidefinite programming relaxations for the graph partitioning problem.
{\em Discrete Appl. Math.}, 96/97:461-–479, 1999.

\end{thebibliography}
\end{document}